\documentclass[preprint,12pt]{elsarticle_clean}

\usepackage[usenames,dvipsnames]{xcolor}

\usepackage{amsmath}
\usepackage{amsfonts}
\usepackage{amssymb}
\usepackage{amsthm}
\usepackage{newlfont}
\usepackage{url}
\usepackage{amssymb}

\newtheorem{theorem}{Theorem}

\newtheorem{definition}{Definition}
 
\newtheorem{assumption}{Assumption}

\newtheorem{remark}{Remark}

\begin{document}
\begin{frontmatter}

\title{Asymptotic behaviour of  dynamical systems with  plastic self-organising vector fields}

\author{N.B. Janson$^{1}$ and P.E.  Kloeden$^{2}$}

\address{$^{1}$Department of Mathematics, Loughborough University, Loughborough LE11 3TU, UK\\
$^{2}$Mathematics Department, University of T\"ubingen, T\"ubingen 72074, Germany}



\begin{abstract}
In [Janson $\&$ Marsden 2017] a dynamical system with a plastic self-organising velocity vector field was introduced, which was inspired by the architectural plasticity of the brain and proposed as a possible conceptual model of a cognitive system. Here we provide a more rigorous mathematical formulation of this problem, make several simplifying assumptions about the form of the model and of the applied stimulus, and perform its mathematical analysis. Namely, we explore the existence, uniqueness, continuity and smoothness of both the plastic velocity vector field controlling the observable behaviour of the system, and of the behaviour itself. We also analyse the existence of pullback attractors and of forward limit sets in such a non-autonomous system of a special form. Our results verify the consistency of the problem, which was only assumed in the previous work, and pave the way to constructing models with more sophisticated cognitive functions. 
\end{abstract}


\end{frontmatter}

\section{Introduction}

It has been widely believed that replication in hardware or software of the brain's physical architecture would automatically replicate the brain's cognitive functions. To a certain extent, this assumption has been correct, since artificial neural networks have been highly successful in a range of applications requiring classification and pattern recognition \cite{Hassabis_brain-inspired_AI_Neu17}. However, they are still far away from demonstrating human-level intelligence. There is a growing appreciation that to recreate the brain's functions, one needs to reveal and reproduce the brain's working principles rather than the architecture per se. To achieve that, it is necessary to answer a number of fundamental questions posed by neuroscience, such as how memories are represented in the brain, how behaviour could be linked to the brain substance, and how cognitive processes could be described in rigorous terms \cite{Damasio_how_brain_creates_mind_SciAm99,Carandini_from_circuits_to_behavior_NN12,Abbott_solving_the_brain_Nat13,Yuste_new_century_of_brain_SA14,Grillner_criticism_neuroscience_Neu14,Katkov_memory_retrieval_Neuron17,Siegfried_long_way_in_understanding_brain_SciNews17}.

Given that the neurons in the brain fire {\it spontaneously}, biologically relevant brain models take the form of dynamical systems with continuous time \cite{Hodgkin_nuron_model_JP52, Leon_virtual_brain_Scholarpedia13}, whose key element is the velocity vector field combining two features. Firstly, assuming that the model of a  spontaneously evolving device is derived from the first principles and is accurate, its velocity field is a mathematical representation of the device's physical architecture  \cite{Janson_conceptual_brain_model_SR17_SN}. Secondly, this field is a mathematical expression of the force fully controlling the device's behaviour, and is an embodiment of the full set of behavioural rules. 

In \cite{Janson_conceptual_brain_model_SR17}  it was proposed that looking at the brain through the prism of its velocity vector field offers a solution to a number of fundamental questions asked by neuroscience. Firstly, the velocity field of the brain could represent the sought-after link between  its physical properties and the resultant behaviour. Secondly, by hypothesising that memories could be imprints on the brain's velocity field, one could unify several dominating memory theories. Thirdly, 
the principles of the brain cognitive function could amount to the ability to create a plastic self-organising velocity vector field evolving according to certain rules, which need to be of an appropriate form to enable cognition. Ultimately, it has been suggested that to be cognitive, the system does not necessarily need to be a neural network, but rather to be capable of spontaneous modifications of its velocity vector field according to some suitable rules. 

Thus, a conceptual model of a cognitive system has been proposed, which represents a dynamical system with a plastic self-organising velocity vector field. The standard part of this model is given by the non-autonomous differential equation:
\begin{equation}\label{xeqn0}
\frac{dx}{dt}    =   a(x,t).
 \end{equation}
Here, $x$ $\in$ $\mathbb{R}^d$ is the state of the cognitive system at any time $t$, which in the brain would be a collection of states of all the neurons, and 
 $a$ $\in$ $\mathbb{R}^d$ is the velocity vector field governing evolution of the state. The unconventional part of the proposed model assumes that the vector field  $a$ evolves with time  
 {\it spontaneously} according to some pre-defined deterministic  rules  and is affected by a stimulus $\eta(t)$ $\in$ $\mathbb{R}^m$, $m$ $\le$ $d$,  as expressed by an equation \begin{equation}\label{aeqn0}
\frac{\partial a}{\partial t}   =  c(a,x,\eta(t),t),  
 \end{equation}
 with $c$ taking values in  $\mathbb{R}^d$.
 
 In \cite{Janson_conceptual_brain_model_SR17} equation \eqref{aeqn0} was interpreted as a (degenerate) partial differential equation (PDE). For the mathematical analysis in this paper it is   more appropriate to regard \eqref{aeqn0} as 
an ordinary differential equation (ODE) for the unknown variable $a$ which depends on the parameter $x$, and to write it as $\frac{da}{dt}=c(a,x,\eta(t),t)$. For an additional clarity, 
\eqref{aeqn0} could be re-written as a parameterised ODE
$$
\frac{d}{dt} a(z,t)  =  c(a(z,t),z,\eta(t),t), \qquad  z \in \mathbb{R}^d, 
$$
since the parameter $z$ does not evolve according to the ODE \eqref{xeqn0}.

Note, that although it has long been  acknowledged that the brain can be regarded as a dynamical system \cite{vanGelder_98}, consideration of the brain within the classical framework of the dynamical systems theory did not fully explain its abilities for cognition and adaptation. 
As a possible reason for  this, it has been suggested that the theory of dynamical systems has not been sufficiently developed and required extensions directly relevant to cognition \cite{Crutchfield_dyn_embod_cogn_BBS98}. 
 In \cite{Janson_conceptual_brain_model_SR17} the velocity vector field with an ability to self-organize has been suggested as an extension of the dynamical systems theory, which could explain adaptation of the behaviour to the environment as a spontaneous modification of the behavioural {\it rules} specified by this field, where this  modification is induced by the self-organised architectural plasticity of the brain.  Under the assumption that the conceptual model \eqref{xeqn0}--\eqref{aeqn0} was mathematically consistent and $a(t)$ stayed smooth for all $t$,  a simple example of \eqref{aeqn0} was constructed, analysed numerically and shown to perform some basic cognition. One could potentially build more sophisticated examples of $c$ leading to  more advanced information-processing functions. 
However, before this can be done, one needs to address the consistency of \eqref{xeqn0}--\eqref{aeqn0}  and the conditions on $c$ and $\eta(t)$ in \eqref{aeqn0}, under which the solution of \eqref{xeqn0} exists and is unique, which has not beed done to date and which is the purpose of the current paper. 

In Section \ref{sec_ps} we formulate a mathematical problem to be solved here.  In Section \ref{sec_exist} we establish the existence and uniqueness of solutions of the system \eqref{Xeqn}--\eqref{Aeqn}. In Section   \ref{sec_asym} we  show that the first equation \eqref{Xeqn}  has a global  non-autonomous attractor. In Section \ref{disc} we discuss the results obtained.

\section{Problem statement}
\label{sec_ps}

We can interpret  \eqref{xeqn0}--\eqref{aeqn0}  as a system of ODEs with an unconventional structure,    
\begin{eqnarray}\label{Xeqn}
\frac{dx(t)}{dt}   & = & a(x(t),t)
\\[1.4ex]
\label{Aeqn}
\frac{d}{dt} a(z,t)  &=&  c(a(z,t),z,\eta(t),t), \qquad  z \in \mathbb{R}^d, 
\end{eqnarray}
with solutions $x(t)$ and  $a(z,t)$ taking values in $\mathbb{R}^d$.   
The solution  $a(t)$  of \eqref{Aeqn} depends on $z$ $\in$ $\mathbb{R}^d$  as a fixed parameter.    
Note that the solution $x(t)$ of the first equation \eqref{Xeqn} is not inserted into  the second equation \eqref{Aeqn},  
i.e. 
equation \eqref{Aeqn}  is decoupled from 
equation \eqref{Xeqn}. Essentially, we need to solve the second   equation \eqref{Aeqn}  first, independently of equation  \eqref{Xeqn},   to obtain the vector field for the first equation \eqref{Xeqn}. 
 
The purpose  of this paper is to provide a more precise mathematical formulation 
and an analysis of the modelling and numerical work in \cite{Janson_conceptual_brain_model_SR17}. 
In particular, we show  that the system \eqref{Xeqn}--\eqref{Aeqn} is well
posed in the sense of the global existence and uniqueness of its solutions. 

It is also important to consider a long-term behaviour of the newly introduced systems, which is usually described by attractors. 
The concept of an attractor has been successfully extended from the autonomous to the standard non-autonomous dynamical systems of the form 
$\frac{dx}{dt}$ $=$ $f(x,t)$, where $f$ is some fixed vector field function  \cite{Kloeden_NDS_book11}. However, the existence of an attractor where the vector field itself evolves {\it spontaneously} according to \eqref{Aeqn} needs to be proved. 
We show that, under a mild dissipativity  assumption,  the non-autonomous  system generated by 
\eqref{Xeqn} 
has a non-autonomous (or random) attractor.

\section{Existence and uniqueness of solutions}
\label{sec_exist}

 In \cite{Janson_conceptual_brain_model_SR17} it was proposed that in 
\eqref{xeqn0}--\eqref{aeqn0} the stimulus $\eta(t)$ is used both to contribute to the modification of the vector field $a$ according to \eqref{aeqn0}, and to regularly reset the initial conditions of \eqref{xeqn0}. Here, we consider a simplified case, in which we allow  $\eta(t)$ only to affect evolution of $a$. Therefore, when considering  system \eqref{Xeqn}--\eqref{Aeqn} 
we will handle equation \eqref{Xeqn} separately from \eqref{Aeqn}, assuming that the vector field $a$ is known. 
 The existence and uniqueness of solutions of equations \eqref{Xeqn} and \eqref{Aeqn}
 require at least a local Lipschitz property of the right-hand sides $a$ and $c$ in the corresponding state variable, while the existence of an attractor in \eqref{Xeqn} requires a dissipativity property. 

Since the equations \eqref{Xeqn}--\eqref{Aeqn} represent a non-autonomous or a random system,  we need to consider  them on the entire time axis $t$ $\in$ $\mathbb{R}$ (see the discussion in Section  \ref{disc} for when this does not hold).  In particular, the vector field  $a$ should be defined for all values of time $t$ $\in$ $\mathbb{R}$. 
Below we formulate our assumption on the stimulus $\eta$. 

\begin{assumption} \label {assumpS1}
$\eta$ $:$ $\mathbb{R}$ $\rightarrow$ $\mathbb{R}^m$  is continuous.
\end{assumption}
This stimulus signal is considered as a given and fixed input in  the model.

\subsection{Existence and uniqueness of the observable behaviour $x(t)$ of \eqref{Xeqn}}

\begin{assumption} \label {assumpA1}
$a$ $:$ $\mathbb{R}^d\times \mathbb{R}$ $\rightarrow$ $\mathbb{R}^d$ and $\nabla_x a(x,t)$ $:$ $\mathbb{R}^d\times \mathbb{R}$ $\rightarrow$ $\mathbb{R}^{d\times d}$  are continuous in both variables $(x,t)$.
\end{assumption}

 This assumption ensures the vector field $a$ is locally Lipschitz in $x$. Hence, by standard theorems (see Walter \cite[Chapter 2]{Walter_ODE_book98}),  there exists a unique solution $x(t)$ $=$ $x(t;t_0,x_0)$ of  the ODE \eqref{Xeqn}  for each initial condition $x(t_0)$ $=$ $x_0$, at least for a short time interval.

\begin{assumption} \label{assumpA2}
  $a$ $:$ $\mathbb{R}^d\times \mathbb{R}$ $\rightarrow$ $\mathbb{R}^d$ satisfies the dissipativity condition  $\left<a(x,t),x\right>$ $\leq$ $-1$  for $\|x\|$ $\geq$ $R^*$ for  some $R^*$. 
\end{assumption}
\noindent (Here $\|a\|$  $=$ $\sqrt{\sum_{i=1}^d a_i^2}$ is the Euclidean norm on $\mathbb{R}^d$ and $<a,b>$  $=$ $\sum_{i=1}^d a_i b_i$ is the corresponding inner product, for vectors $a$, $b$ $\in$  $\mathbb{R}^d$.)

 This assumption (which may be stronger than we really need, but avoids assumptions about the specific structure of $a$) ensures that the ball $B^*$ $:=$   $\{x \in \mathbb{R}^d :  \|x\|  \leq  R^*+1\}$  is positive invariant. This follows  from the estimate 
$$
\frac{d}{dt}  \|x(t)\|^2 = 2 \left<x(t),a(x(t),t) \right>  \, \,\leq \,\, -1 \qquad \mbox{if}\,\, \|x(t)\| \geq R^*
$$
and in turn ensures that the solution of the ODE \eqref{Xeqn} exists for all future time $t$  $\geq$ $t_0$.
We thus formulate the following theorem. 
\begin{theorem} \label{thm1}
Suppose that Assumptions \eqref{assumpS1}, \eqref{assumpA1} and \eqref{assumpA2} hold. Then for every initial condition $x(t_0)$ $=$ $x_0$, the  ODE \eqref{Xeqn} has a  unique solution
 $x(t)$ $=$ $x(t;t_0,x_0)$, which exists for all $t$ $\geq$ $t_0$.  Moreover, these solutions are  continuous in the initial conditions, i.e., the mapping  $(t_0,x_0)$ $\mapsto$ 
$x(t;t_0,x_0)$ is continuous.
 \end{theorem}

\subsection{Existence and uniqueness of the vector field $a(x,t)$ as a solution of \eqref{Aeqn}}

The ODE   \eqref{Aeqn}  for the 
 velocity field  $a(x,t)$  is independent of the 
 solution $x(t;t_0,x_0)$  of  the ODE \eqref{Xeqn}. We need the following assumption  to provide the existence and uniqueness of   $a(x,t)$   for all future times $t$ $>$ $t_0$ and to ensure  that  this solution satisfies  Assumptions  \eqref{assumpA1} and \eqref{assumpA2}.

\begin{assumption} \label {assumpC1}
$c$ $:$ $\mathbb{R}^d\times \mathbb{R}^d \times \mathbb{R}^m\times \mathbb{R}$ $\rightarrow$ $\mathbb{R}^d$ and $\nabla_a c $ $:$ $\mathbb{R}^d  \times \mathbb{R}^d \times \mathbb{R}^m\times \mathbb{R}$  $\rightarrow$ $\mathbb{R}^{d\times d}$  are continuous in all variables.
\end{assumption}

 This assumption ensures the vector field $c$ is locally Lipschitz in $a$. Hence, by standard theorems (see Walter \cite[Chapter 2]{Walter_ODE_book98}),  there exists a unique solution $a(t;t_0,a_0)$ of  the ODE \eqref{Aeqn}  for each initial condition $a(t_0)$ $=$ $a_0$, at least for a short time interval.  This solution also  depends  continuously on the parameter $x$ $\in$ $\mathbb{R}^d$. 
To ensure that the solutions can be extended for all future times $t$, we need a growth bound such as in the following assumption. 
\begin{assumption} \label {assumpC2}
 There exist constants $\alpha$ and $\beta$ (which need not be positive) such that $\left<a, c(a,x,y,t)\right>$ $\leq$ $\alpha \|a\|^2 +\beta$ for all   $(x,y,t)$ $\in$ $\mathbb{R}^d \times \mathbb{R}^m\times \mathbb{R}$.
\end{assumption}

The next assumption ensures that the solution of the ODE \eqref{Aeqn}, which we now write as $a(x,t)$,   is continuously differentiable and hence locally Lipschitz in $x$, provided  
that  the  initial value $a(x,t_0)$ $=$ $a_0(x)$ is continuously differentiable. 
\begin{assumption} \label {assumpC3}
 $\nabla_x c $ $:$ $\mathbb{R}^d \times \mathbb{R}^d \times \mathbb{R}^m\times \mathbb{R}$  $\rightarrow$ $\mathbb{R}^{d\times d}$ is  continuous in all variables.
\end{assumption}
The above statement then follows  from the properties of the linear matrix-valued  variational equation
$$
\frac{d}{dt}  \nabla_x a =  \nabla_a c   \nabla_x  a +   \nabla_x c,  
$$  
which is obtained by taking the gradient $\nabla_x $ of both sides of the ODE  \eqref{Aeqn}.

Finally, we need to ensure that the solution $a(x,t)$  satisfies the dissipativity  property as in Assumption \eqref{assumpA2}. 

\begin{assumption} \label{assumpC4}
  There exist $R^*$   such that 
  $$
\left<c(a,x,y,t),x\right>  \,  \leq \, 0 \quad  \textrm{for} \quad \|x\|   \geq R^*,   \quad (a,y,t)    \in    \mathbb{R}^d \times \mathbb{R}^m\times \mathbb{R}.
$$
   \end{assumption}
   
To show this we write equation  \eqref{Aeqn} in integral form
$$
a(x,t) = a_0(x) + \int_{t_0}^t c\left(a(x,s),x, \eta(s),s\right)\, ds
$$
and then take the scalar product on  both sides with a constant $x$, which gives
\begin{eqnarray*}
\left<a(x,t),x\right> & = &  \left<a_0(x),x\right> + \left< \int_{t_0}^t   c\left(a(x,s),x, \eta(s),s\right)  \, ds, x\right>
\\[1.3ex] 
& = &  \left<a_0(x),x\right> + \int_{t_0}^t \left< c\left(a(x,s),x, \eta(s),s\right),x\right>\, ds
\\[1.3ex] 
& \leq &   - 1 +0 = -1 \quad  \mbox{\rm  for} \qquad \|x\|   \geq R^*.
\end{eqnarray*}

Summarising from the above, we can formulate the following theorem. 

\begin{theorem}
 \label{thm2}
Suppose that Assumptions \eqref{assumpS1} and \eqref{assumpC1}--\eqref{assumpC4}  hold.  Further, suppose that  $a_0(x)$  is continuously differentiable and satisfies the dissipativity condition in Assumption  \eqref{assumpA2}.  Then the ODE \eqref{Aeqn} has a unique  solution $a(x,t)$ for the initial condition $a(x,t_0)$ $=$ $a_0(x)$, 
which exists  for all $t$ $\geq$ $t_0$ and satisfies   Assumptions  \eqref{assumpA1} and \eqref{assumpA2}. 
 \end{theorem}
Thus, we have obtained a theorem for the existence, uniqueness, continuity and dissipativity of the velocity vector field $a$ governing the behaviour of \eqref{Xeqn}.

\section{Asymptotic behaviour}

\label{sec_asym}
 
Here we consider the conditions for the existence of two kinds of attractors in equation \eqref{Xeqn} describing the observable behaviour of a system with a plastic velocity field. 
The  ODE \eqref{Xeqn}   is non-autonomous and its solution mapping   generates a non-autonomous dynamical system on the state space $\mathbb{R}^d$ expressed in terms of a $2$-parameter  semi-group, which is often called a  process (see Kloeden \& Rasmussen \cite{Kloeden_NDS_book11}).  Define 
\begin{equation*}
  \mathbb{R}_{\geq}^+=
  \left\{(t,t_0)\in\mathbb{R}\times\mathbb{R}:t\geq t_0\right\}. 
\end{equation*}

\begin{definition}  \label{DETdefpro}
  A \emph{process} is a mapping $\phi:\mathbb{R}_{\geq}^+\times \mathbb{R}^d\to \mathbb{R}^d$ with the 
  following properties: 
  \begin{itemize} 
  \item[(i)]   {initial condition:}\,\, 
    $\phi(t_0,t_0,x_0)=x_0$ for all $x_0\in \mathbb{R}^d$ and $t_0\in\mathbb{R}$;
  \item[(ii)]  {$2$-parameter semi-group property:} \, 
    $\phi(t_2,t_0,x_0)=\phi(t_2,t_1,\phi(t_1,t_0,x_0))$ 
    for all    $t_0\leq t_1  \leq t_2$ in $\mathbb{R}$ and $x_0\in \mathbb{R}^d$;
  \item[(iii)]   {continuity:}\,\, 
    the mapping $(t,t_0,x_0)\mapsto\phi(t,t_0,x_0)$ is continuous.
  \end{itemize}
\end{definition}
The $2$-parameter semi-group property is an immediate consequence  of the existence and uniqueness of solutions of the non-autonomous 
ODE: the solution starting  at~$(t_1,x_1)$, where $x_1$ $=$ $\phi(t_1,t_0,x_0)$, is unique so must be equal to $\phi(t,t_0,x_0)$ for  $t$ $\geq$ $t_1$.

\subsection{Pullback attractors in equation \eqref{Xeqn}} 

Time in an autonomous dynamical systems is a relative concept since such systems depend   on the elapsed time $t-t_0$ only and not separately on the 
current time~$t$ and initial time~$t_0$, which means that limiting  objects exist  all the time and not just in the distant future.  In contrast, non-autonomous systems depend explicitly on both~$t$  and~$t_0$, which has a profound affect on the nature of limiting objects (see \cite{Kloeden_NDS_book11,
Crauel_nonaut_random_attr2015}).

In particular,   a non-autonomous   attractor is a family   $\mathfrak{A}$ $=$ $\{ A(t): t \in \mathbb{R}\}$ of  nonempty compact subsets  $A(t)$ of $\mathbb{R}^d$  with  the following properties:\\ 

\noindent 1) invariance: $A(t)$ $=$ $\phi(t,t_0,A(t_0))$  for all $t$ $\geq$ $t_0$;  

\noindent 2) pullback attracting: 
$$
\lim_{t_0\to - \infty} \mbox{\rm dist}_{\mathbb{R}^d} \left(\phi(t,t_0,B),A(t)\right) = 0  \qquad \mbox{for all bounded subsets $B$ of    $\mathbb{R}^d$. }\\
$$
It is   called  a \emph{pullback attractor} since  the starting time $t_0$ is pulled further and further back into the past. The dynamics then moves forwards in time from this starting time $t_0$ 
to the present time $t$. Essentially, the pullback attractor takes into account the past history of the system, so we cannot expect it to say much about the future. In fact, a pullback  attractor need not be forward attracting in the  conventionally understood sense, i.e. as $t$ $\to$ $+\infty$ for fixed $t_0$  (see the example in subsection  in \ref{subsec_example} below).

The existence and uniqueness of a global pullback attractor for a non-autonomous dynamical system on $\mathbb{R}^d$ is  implied by the existence of a positive invariant absorbing set. 
The following  theorem is adapted from \cite[Theorem  3.18]{Kloeden_NDS_book11}. 
\begin{theorem}\label{thm3}
Suppose that a   non-autonomous dynamical system $\phi$ on $\mathbb{R}^d$ has  a positive invariant absorbing set $B^*$. Then it has a unique pullback attractor  $\mathfrak{A}$ $=$ $\{ A(t): t \in \mathbb{R}\}$  with component sets 
defined by
$$
A(t)  = \bigcap_{t_0\leq t}  \phi\left(t,t_0,B^*\right), \qquad  t \in \mathbb{R}.  
$$
\end{theorem}
An important characterization \cite[Lemma 2.15]{Kloeden_NDS_book11} of a pullback attractor is that it consists of the entire bounded solutions of the system, i.e., $\chi$ $:$ $\mathbb{R}$ $\rightarrow$ $\mathbb{R}^d$ for which $\chi(t)$ $=$ $\phi(t,t_0,\chi(t_0))$ $\in$ $A(t)$ for all $(t,t_0)$ $\in$   $\mathbb{R}_{\geq}^+$.

In particular, under the above assumptions, the ODE  \eqref{Xeqn}  describing the observable  behaviour of the model of a cognitive system generates a  non-autonomous dynamical system, which has a global pullback attractor.
Summarising, we formulate the following theorem. 
\begin{theorem} \label{thm4}
Suppose that Assumptions \eqref{assumpS1}, \eqref{assumpA1} and \eqref{assumpA2} hold. Then  the  non-autonomous dynamical system generated by the ODE \eqref{Xeqn} describing the observable behaviour  
has a   global pullback attractor  $\mathfrak{A}$ $=$ $\{ A(t): t \in \mathbb{R}\}$, which is   contained in the absorbing set $B^*$.
 \end{theorem}
 Thus, Theorem \ref{thm4} specifies the conditions under which the global pullback attractor exists in a dynamical system with plastic spontaneously evolving velocity vector field. 
 
\subsection{Forward limit sets in equation \eqref{Xeqn}}  
The concepts of pullback attraction and pullback attractors assume that the system  exists for all time, in particular past time.  This is obviously not true  in many biological systems, though 
an artificial ``past" can some times be usefully introduced (see the final section).  

The  above definition of a non-autonomous dynamical system can be easily modified to hold only for $(t,t_0)$ $\in$    $\mathbb{R}_{\geq}^+(T^*)$ $=$ $\left\{(t,t_0)\in\mathbb{R}\times\mathbb{R}: t\geq t_0\geq T^* \right\}$  for some $T^*$ $>$ $-\infty$. 

When the system has  a nonempty positive invariant compact   absorbing set $B^*$, as in  the situation here,  the forward omega limit set   
$$
\Omega(t_0) = \bigcap_{\tau \geq t_0} \overline{ \bigcup_{t\geq \tau}  \phi\left(t,t_0,B^*\right)}, \qquad  t_0 \in \mathbb{R}, 
$$
 exists for each $t_0$ $\geq$ $T^*$, where the upper bar denotes the closure of the set under it.  The set  $\Omega(t_0)$ is thus a nonempty 
 compact subset  of the absorbing set $B^*$ for each $t_0$ $\in$ $\mathbb{R}$.
 
 Moreover,   these sets are increasing in $t_0$, i.e., $\Omega(t_0)$ $\subset$ $\Omega(t'_0)$ for $t_0$ $\leq$ $t'_0$, and the closure of their union 
 $$
 \Omega^* :=   \overline{ \bigcup_{t_0\geq T^*}   \Omega(t_0)}  \subset B^*
$$
is a compact subset of $B^*$, which attracts all of the dynamics of the system in the forward sense, i.e.  
$$
\lim_{t\to   \infty} \mbox{\rm dist}_{\mathbb{R}^d} \left(\phi(t,t_0,B), \Omega^* \right) = 0  
$$ 
for all bounded subsets $B$ of  $\mathbb{R}^d$, $t_0$ $\geq$ $T^*$. 
 
Vishik \cite{Vishik_Asym_Beh_book92}     called    $\Omega^*$ the \emph{uniform attractor}\footnote{He required the system to be defined   in the whole past and the convergence to be uniform in $t_0$ $\in$ $\mathbb{R}$.},  although strictly speaking   $\Omega^*$  do not form an attractor  since it need not be invariant and the attraction need not be uniform in the starting time $t_0$.   Nevertheless,   $\Omega^*$   does indicate where the future asymptotic dynamics ends up. Moreover,  Kloeden  \cite{Kloeden_asym_inv_JCD16} showed that $\Omega^*$  is  \emph{asymptotically positive invariant},  which means that  the later the starting time $t_0$, the more and more it looks like an attractor as conventionally understood.   In  \cite{Kloeden_asym_inv_JCD16} $\Omega^*$  was called the  \emph{forward attracting set}.

Summarising   from the above, we formulate the following theorem. 
\begin{theorem} \label{thm5}
Suppose that Assumptions \eqref{assumpS1}, \eqref{assumpA1} and \eqref{assumpA2} hold. Then  the  non-autonomous dynamical system generated by the  
ODE \eqref{Xeqn} describing the observable behaviour of the system 
has a   forward attracting set   $\Omega^*$, which is  contained in the absorbing set $B^*$.
 \end{theorem}
 
 Theorem \ref{thm5} expresses the conditions under which a forward attracting set exists in a dynamical system \eqref{Xeqn}  with a plastic  velocity vector field evolving according to  \eqref{Aeqn}.

 \section{Discussion} \label{disc}
 
In the previous sections we explicated the mathematical formulation of, and analysed mathematically, the conceptual model of a cognitive system introduced in   \cite{Janson_conceptual_brain_model_SR17}. For clarity, we made some simplifying assumptions about the properties of the right-hand sides of this model and of the external stimulus. 

If the model discussed here is to be used for the description of the cognitive function similar to that of a biological brain, one needs to take into account the different timescales at which different processes occur. It is known that the observable dynamics of neurons is much faster than the rate of change of the inter-neuron connections. 
Hence, the velocity vector field describing the dynamics of neurons in the brain should evolve at a much slower rate than the neural states. A realistic application of our model should take this into account. 

Even the simplified cases studied here raise a number of questions, in particular  about the relevance of pullback attractors for such models. These and some further issues will be briefly discussed here.

  \subsection{Use of pullback attractors}

 Pullback convergence requires the dynamical system to exist in the  distant past, which is often not a realistic assumption  in biological systems. Pullback attractors can nevertheless be used in such situations by inventing an artificial past. This  and other  aspects are discussed in  \cite{Kloeden_pullback_SD03,Kloeden_NDS_book11}. 
 
 The simplest way to do this  for this model is to set the  vector field $a(x,t)$ $\equiv$  $a_0(x)$    for $t$ $\leq$ $T^*$   for some   finite time $T^*$, which could be the 
 desired  starting time $t_0$. In this case $a_0(x)$ would be the desired initial velocity vector field of the model of a cognitive system, which could be zero or contain some initial features representing  previous memories. Then the ODE \eqref{xeqn0} should be replaced by the  switching system
 \begin{eqnarray}
 \label{eq_sw}
\frac{dx}{dt} = \left\{ \begin{array}{lcl} a_0(x)   & : & t \leq t_0   \\[1.3ex]   a(x,t) & : & t \geq t_0  
\end{array}\right.,
\end{eqnarray}
where $a(x,t)$ evolves  according to the ODE \eqref{aeqn0} { for $t$ $\geq$ $t_0$  with the parameterised initial value $a_0(x)$. If $a_0(x)$ satisfies the dissipativity condition in Assumption  \eqref{assumpA2}, then the switching system \eqref{eq_sw} will also be dissipative and have a pullback attractor with component  sets $A(t)$ $=$ $A^*$ for $t$ $\leq$ $t_0$ and $A(t)$ $=$ 
$\phi(t,t_0,A^*)$ for $t$ $\geq$ $t_0$, where $A^*$ is the global attractor of the autonomous dynamical system generated by the autonomous ODE with the vector field   $a_0(x)$.  

\smallskip

\subsection{Random stimulus signals}
  
The stimulus signal $\eta(t)$  in Assumption \eqref {assumpS1} is a deterministic  function. When this signal is random it would be a single sample path  $\eta(t,\omega)$ of a stochastic process  with  $\omega$ $\in$ $\Omega$, where $\Omega$ is   the sample space of the underlying probability space$(\Omega, \mathcal{F}, \mathbb{P})$.  The above analysis holds, which is otherwise deterministic,  for this  fixed sample path. For emphasis,  $\omega$ could be included in the system and the pullback attractor as an additional parameter, i.e.,  $\phi(t,t_0,x_0, \omega)$ and 
$\mathfrak{A}$ $=$ $\{ A(t,\omega): t \in \mathbb{R}\}$. Cui et al \cite{Cui_forward_random_att_submitted,Cui_uniform_att_nonaut_JDE17} call these objects non-autonomous random dynamical systems and random pullback attractors,  respectively. 

This is the appropriate formulation for brain  vector fields generated by the ODE  \eqref{Aeqn}, which has two sources of non-autonomity  in its vector field $c$, i.e., indirectly through the stimulus signal  $\eta(t)$ and directly through the  independent variable $t$. The  ODE \eqref{Aeqn} is then a random ODE (RODE), see \cite{Han_RODE_book17}. Note, that without this additional  independent variable $t$, the  theory of random dynamical  systems (RDS)  in Arnold \cite{Arnold_RDS_book98} could be used. It is also a pathwise theory with a random attractor defined through pullback convergence, but  requires additional assumptions about the nature of the driving noise process, which is here represented in the stimulus signal.

Until now we considered the stimulus signal $\eta(t)$ with continuous sample paths. The above results remain valid when $\eta(t)$ has only  measurable sample  paths,  such as for a Poisson or L\'evy process,   but the RODE  must now be interpreted pathwise as  a Carath\'eodory ODE,  see \cite{Han_RODE_book17}.

 \subsection{Relevance of pullback attractors}
 \label{subsec_example}
   
Assuming,  by nature or artifice, that the system does have a pullback attractor, what does this  actually tell us about the asymptotic dynamics of the neural activity?
  
As mentioned above,  a pullback attractor consists of the entire bounded solutions of the system, which is useful information.   This characterisation is also true of attractors of autonomous systems, for which  pullback and forward convergence are equivalent due to the fact that only the elapsed time is important in such systems.  

In general,  a pullback attractor need not be forward attracting. This is easily seen in the following switching system
$$
\frac{dx}{dt}=   \left\{ \begin{array}{ccl} - x  & : & t \leq 0 
\\[1.5ex]
x \left(1-x^2 \right) &   : & t > 0 
\end{array} \right. ,
$$
for which the set $B^*$ $=$ $[-2,2]$ is  positively invariant and absorbing. The  pullback attractor   $\mathcal{A}$   has identical component subsets  $A_t$ $\equiv$ $\{0\}$, $t$ $\in$ $\mathbb{R}$,  corresponding to the zero entire solution, which is  the only bounded entire solution of this switching  system.  This zero solution is  obviously  not forward asymptotically stable. The  forward attracting set here is $\Omega^*$ $=$ $[-1,1]$. It is not invariant (though it is  positive invariant  in  this case), but  contains all of the forward limit points of the system.
 
Nevertheless,  a pullback attractor   indicates where the system settles down to  when more and more information of its past is taken into account. This is very  useful in system which is itself  evolving in time,  as in   the brain plasticity model under consideration, for which the future input stimulus is not yet known. 

 Interestingly, a random attractor for the  RDS  in  the sense of \cite{Arnold_RDS_book98} is  pullback attracting in the pathwise sense and also forward attracting in probability, see \cite{Crauel_nonaut_random_attr2015}. 
 
\subsection{Vector field  from a potential function}

In an example investigated numerically in  \cite{Janson_conceptual_brain_model_SR17},  the vector field $a$ was generated  from  a potential function $U$, i.e.  with $a$ $=$ 
$-\frac{1}{t}\nabla_x U$, and a differential equation was constructed  for  $U$,  rather than for $a$. Componentwise $a_i$ $=$ $-\frac{1}{t}\frac{\partial U}{\partial x_i}$, so the existence of such a potential requires  
$$
\frac{\partial a_i}{\partial x_j} = -\frac{1}{t}\frac{\partial^2 U}{\partial x_j\partial x_i} = -\frac{1}{t}  \frac{\partial^2 U}{\partial x_i\partial x_j} = \frac{\partial a_j}{\partial x_i}
$$
From equation \eqref{Aeqn} this requires 
$$
 \frac{\partial a_i}{\partial x_j}  = \frac{\partial a_j}{\partial x_i}.
$$

The example considered in \cite{Janson_conceptual_brain_model_SR17} is a special case of the system \eqref{Xeqn}--\eqref{Aeqn}. 
Namely, 
 in \cite{Janson_conceptual_brain_model_SR17} $U$ satisfies a scalar parameterised  ordinary differential equation
\begin{equation}\label{ueq}
\frac{d}{dt} U(x,t)  = - k U(x,t)  -g(x-\eta(t)),
\end{equation}
where $k$ $\geq$ $0$,  $g$ is shaped like a  Gaussian function 
$$
g(z) =  \frac{1}{\sqrt{2 \pi \sigma^2}} e^{-\frac{z^2}{ \sigma^2} },  
$$
and $\eta(t)$ is the given input, which is assumed to  be defined for all $t \in \mathbb{R}$ and is continuous.

The gradient $\nabla_x U$ of $U$ satisfies  the scalar   parameterised  ordinary differential equation 
\begin{equation}\label{ueq_grad}
\frac{d}{dt} \nabla_x U(x,t)  =  - k  \nabla_x U(x,t)  - G(x-\eta(t)),
\end{equation}
where 
$$
G(x-\eta(t)) =   \nabla_x  g(x-\eta(t)) = - \frac{2}{\sigma^2\sqrt{2 \pi \sigma^2}} \left(x-\eta(t)\right) e^{-\frac{\left(x-\eta(t)\right) ^2}{ \sigma^2} }.
$$
The  linear ODE \eqref{ueq_grad} has an explicit solution 
$$
 \nabla_x U(x,t) =   \nabla_x U(x,t_0) e^{-k(t-t_0)} - \int_{t_0}^t e^{-k(t-s)} G(x-\eta(s)) ds.
$$
Taking the pullback limit as $t_0$ $\to$ $-\infty$  gives
$$
 \nabla_x \bar{U}(x,t) =     - \int_{-\infty}^t e^{-k(t-s)} G(x-\eta(s)) ds.
$$ 
This solution is asymptotically stable and forward attracts all other solutions,  since 
$$
\left|  \nabla_x U(x,t) -  \nabla_x \bar{U}(x,t)  \right| 	 \leq \left|  \nabla_x U(x,t_0) -  \nabla_x \bar{U}(x,t_0)  \right| e^{-k(t-t_0)}
$$
for every $x$ and any solution $ \nabla_x U(x,t)$ $\ne$  $\nabla_x \bar{U}(x,t)$.

Finally, the asymptotic dynamics of this example system with a plastic vector field satisfies the scalar ODE
\begin{equation}
\label{sol_asym}
\frac{dx(t)}{dt}  = - \frac{1}{t} \nabla_x \bar{U}(x(t),t) =     \frac{1}{t} \int_{-\infty}^t e^{-k(t-s)} G(x(t)-\eta(s)) ds.
\end{equation}
Since the integrand is uniformly bounded, it follows that $\left|\frac{dx(t)}{dt}\right|$ $\leq$  $\frac{C}{t}$ 
$\to$ $0$ as  $t$ $\to$ $\infty$. From numerical simulations, the system \eqref{Xeqn} with $a(x,t)$ $=$ 
$-\frac{1}{t}\nabla_x U(x,t)$ appears to have  a forward attracting set.

From the argument presented above, equation \eqref{Aeqn}  for the vector field $a$ has a pullback attractor consisting of singleton set, i.e. a single entire solutions, which is also Lyapunov forward attracting. 
This implies that starting from an arbitrary smooth initial vector field $a(x,t_0)$, the solution $a(x,t)$ of  \eqref{Aeqn} converges to 
 a time-varying  function $\bar{a}(x,t)$ $=$ $ - \frac{1}{t} \nabla_x \bar{U}(x,t)$.

\begin{remark} The example considered in \cite{Janson_conceptual_brain_model_SR17} actually involved a random forcing term $\eta(t)$, which was  the 
stochastic  stationary solution (essentially its random attractor) of the scalar It\^o stochastic differential equation (SDE)
\begin{equation}\label{noise}
d\eta (t) = h(\eta(t)) dt + 0.5 dW(t),
\end{equation}
where $W(t)$ was a two-sided Wiener process. For the function $h(u)$ $=$  $3(u-u^3)/5$ used in \cite{Janson_conceptual_brain_model_SR17}, the representative potential function had two non-symmetric wells of different depths and widths.   In such cases, the solutions of \eqref{Xeqn}--\eqref{Aeqn} 
depend on the sample path $\eta(t,\omega)$  of the noise  process, and the convergences are pathwise, and  random versions of the   theorems formulated above  
apply.   In particular, the random pullback attractor consists of singleton sets, i.e., it is essentially is a stochastic process. Moreover, it is 
 Lyapunov asymptotically stable in probability.
\end{remark}

\section{Conclusion}

To conclude, we reconsidered the problem from \cite{Janson_conceptual_brain_model_SR17} from the perspective of recent developments in non-autonomous dynamical systems. In order to further develop modelling of information processing by means of dynamical systems with plastic self-organising vector fields, we needed to show that the problem is well-posed mathematically, which is one of the results of this paper obtained under some simplifying assumptions. 
At the same time, we have shown that asymptotic dynamics can be formulated in terms of non-autonomous and/or random attractors. This provides us with a firm foundation for a deeper understanding of the potential capabilities of systems with plastic adaptable rules of behaviour. 

The model presented here offers many interesting mathematical challenges, such as the rigorous analysis of parameteter-free bifurcations occurring as a result of spontaneous evolution of the velocity field of the dynamical system. The necessary background theory is yet to be developed.  

~\vspace{-1mm}

\noindent {\bf Acknowledgements}

\noindent The visit of PEK to Loughborough University was supported by London Mathematical Society. 

~\vspace{-1mm}

\noindent {\bf  Authors' Contributions}

\noindent  NBJ and PEK jointly formulated the problem and interpreted the results. PEK obtained the mathematical results and NBJ put the problem into the context of applications to the brain and cognition. 








\bibliographystyle{elsarticle-num}
\bibliography{brain_model}

\end{document}